\documentclass[11pt,twoside]{article}

\usepackage{amsmath,amssymb,amsthm,url}

\newcommand{\Prim}{\mathrm{Prim}}

\newcommand{\vol}{\mathrm{vol}}

\newcommand{\Tr}{\mathrm{Tr}}
\newcommand{\tr}{\mathrm{tr}}

\newcommand{\Ind}{\mathrm{Ind}}

\newcommand{\as}{\quad\text{as}\quad}
\newcommand{\tinf}{\to\infty}
\newcommand{\disp}{\displaystyle}
\newcommand{\bsla}{\backslash}

\newcommand{\cO}{\mathcal{O}}

\newcommand{\bC}{\mathbb{C}}
\newcommand{\bR}{\mathbb{R}}
\newcommand{\bQ}{\mathbb{Q}}
\newcommand{\bZ}{\mathbb{Z}}

\newcommand{\noi}{\noindent}
\renewcommand{\Re}{\mathrm{Re}}
\renewcommand{\Im}{\mathrm{Im}}

\newcommand{\divset}{\hspace{3pt}|\hspace{3pt}}

\newcommand{\gam}{\gamma}
\newcommand{\Gam}{\Gamma}

\newcommand{\bLambda}{\bar{\Lambda}}

\newcommand{\slt}{\mathrm{SL}_2}

\newcommand{\sr}{\mathrm{SL}_2(\bR)}
\newcommand{\sz}{\mathrm{SL}_2(\bZ)}

\newcommand{\vcpt}{\vspace{12pt}}

\newtheorem{thm}{Theorem}[section]
\newtheorem{prop}[thm]{Proposition}
\newtheorem{lem}[thm]{Lemma}

\numberwithin{equation}{section}

\title{Square integrals of the logarithmic derivatives of 
Selberg's zeta functions in the critical strip}
\author{Yasufumi Hashimoto}
\date{}

\begin{document}
\markboth
{Y. Hashimoto}
{Square integrals of Selberg's zeta functions}
\pagestyle{myheadings}

\maketitle
\renewcommand{\thefootnote}{}

\begin{abstract} 
In our previous work (\url{https://doi.org/10.1002/mana.202000268}, Math. Nachr., 2021), 
we proposed an upper bound of the logarithmic derivative of Selberg's zeta function 
for the modular groups in the critical strip. 
The present paper studies the growth of its square integral for the modular group, 
co-compact arithmetic groups derived from indefinite quaternion algebras 
and their subgroups.
\end{abstract}

\section{Introduction and the main theorem}
Let $H:=\{x+y\sqrt{-1}\divset x,y\in\bR,y>0\}$ be the upper half plane 
and $\Gam$ a discrete subgroup of $\sr$ with $\vol(\Gam\bsla H)<\infty$. 
Denote by $\Prim(\Gam)$ the set of primitive hyperbolic conjugacy classes 
of $\Gam$ and $N(\gam)$ the square of the larger eigenvalue of $\gam$. 
The Selberg zeta function for $\Gam$ is defined by 
\begin{align*}
Z_{\Gam}(s):=\prod_{\gam\in\Prim(\Gam),n\geq0}
(1-N(\gam)^{-s-n}),\qquad \Re{s}>1.
\end{align*}
It is well known that $Z_{\Gam}(s)$ is analytically continued to the whole complex plane 
as a meromorphic function 
and has a functional equation between the values at $s$ and $1-s$ (see, e.g. \cite{He}). 
In the present paper, we study the growth of $Z_{\Gam}(s)$ 
in the critical strip $\{0<\Re{s}<1\}$ as $|\Im{s}|\tinf$.

For the Riemann zeta function 
$\zeta(s):=\prod_{p}(1-p^{-s})^{-1},$ $(\Re{s}>1),$
it has been known there exists a constant $a\geq0$ satisfying 
$$\zeta\left(\frac{1}{2}+iT\right)\ll_{\epsilon}T^{a+\epsilon}, \as T\tinf.$$ 
Note that  $a\leq \frac{32}{205}$ was proven \cite{Hux} 
and $a=0$ has been expected as the Lindel\"{o}f conjecture. 
For the Selberg zeta function, 
the following bounds have been given for for $1/2<\sigma<1$ (see e.g. \cite{He,IwPGT,AS}).
\begin{align}
\log{Z_{\Gam}(\sigma+iT)},\quad  \frac{Z'_{\Gam}(\sigma+iT)}{Z_{\Gam}(\sigma+iT)}
\ll_{\epsilon}T^{2-2\sigma+\epsilon}, 
\as T\tinf.  \label{SZorg}
\end{align}
The bound above means that  
the growth of $Z_{\Gam}(s)$ is (presently bounded to be) exponential of $|\Im{s}|$,  
which is quite different to the Riemann zeta function.
This fact often causes difficulties 
when analyzing the Selberg zeta functions, 
and then sharpening the estimate \eqref{SZorg} has been desirable. 

In our previous work \cite{HashSZ}, 
we improved \eqref{SZorg} for $\Gam=\sz$ to 
\begin{align}\label{hash}
\frac{Z'_{\Gam}(\sigma+iT)}{Z_{\Gam}(\sigma+iT)}\ll_{\epsilon}
\begin{cases}
\disp T^{\frac{19}{9}-\frac{20}{9}\sigma+\epsilon}, & (\frac{1}{2}<\sigma\leq \frac{5}{8}),\\
\disp T^{\frac{52}{27}(1-\sigma)+\epsilon}, & (\frac{5}{8}< \sigma<1),
\end{cases} \as T\tinf
\end{align}
by using the expression of the Selberg zeta function for $\sz$ 
in terms of indefinite binary quadratic forms. 
On the other hand, Aoki \cite{Aoki} studied the square integral 
$$\int_{1}^{T}\left|\frac{Z'_{\Gam}(\sigma+it)}{Z_{\Gam}(\sigma+it)}\right|^{2}dt$$
and proved that it appears in the pair correlation formulas 
for the non-trivial zeros of the Selberg zeta functions 
under the assumption that 
$\Gam$ is co-compact torsion free 
and $Z_{\Gam}(s)$ does not have real zeros in the critical strip. 
He also proved that there exists a constant $C_{\Gam}(\sigma)>0$ such that 
\begin{align}\label{ao}
\int_{1}^{T}\left|\frac{Z'_{\Gam}(\sigma+it)}{Z_{\Gam}(\sigma+it)}\right|^{2}dt
\sim C_{\Gam}(\sigma)T, \as T\tinf 
\end{align}
when $\frac{17}{20}<\sigma <1$ and 
$\Gam$ is a co-compact arithmetic group derived from an indefinite division quaternion algebra, 
defined as follows: 
Let $a,b$ be square free and relatively prime positive integers and  
$B:=\bQ+\bQ\alpha+\bQ\beta+\bQ\alpha\beta$ the quaternion algebra over $\bQ$
with $\alpha^2=a$, $\beta^2=b$, $\alpha\beta=-\beta\alpha$. 
Suppose that $B$ is division and  
fix a maximal order $\mathcal{O}$ of $B$. 
Then the group $\mathcal{O}^1$ consisting 
of $q_0+q_1\alpha+q_2\beta+q_3\alpha\beta\in \cO$ 
with $q_0^2-q_1^2a-q_2^2b+q_3^2ab=1$  
can be identified with 
a co-compact discrete subgroup $\Gam=\Gamma_{\mathcal{O}}$ 
of $\sr$ by the map
\begin{align*}
q_0+q_1\alpha+q_2\beta+q_3\alpha\beta\mapsto &
\begin{pmatrix} 
q_0+q_1\sqrt{a}&q_2\sqrt{b}+q_3\sqrt{ab}\\
q_2\sqrt{b}-q_3\sqrt{ab}&q_0-q_1\sqrt{a}
\end{pmatrix}.
\end{align*}

While \eqref{hash} and \eqref{ao} are for different $\Gam$'s, 
the asymptotic formula \eqref{ao} gives a better estimate
than the square integral of the right hand side of \eqref{hash} 
in the region $\{\frac{17}{20}<\sigma < 1\}$. 
It was proven by 
Landau's formula for square integrals of Dirichlet series 
(see Section 233--226 in \cite{Lan}) based on the fact that 
the error term of the prime geodesic theorem for the corresponding quaternion $\Gam$ 
are bounded as follows (see \cite{Ko,LRS}).
\begin{align*}
\#\{\gam\in\Prim(\Gam) \divset N(\gam)<x\}-\int_{2}^{x}\frac{dt}{\log{t}}
\ll_{\epsilon,\Gam} x^{\frac{7}{10}+\epsilon},\as x\tinf.
\end{align*} 
The range $\frac{17}{20}<\sigma < 1$ can be extended 
if the error term estimation of the prime geodesic theorem will be improved.
For the modular group, the exponent $\frac{7}{10}$ of the error term 
was improved to $\frac{25}{36}$, and to $\frac{2}{3}$ 
if the Lindel\"{o}f conjecture 
for the Dirichlet $L$ function holds \cite{SY,BF}. 
By using Landau's formula, one can improve the range of $\sigma$ for \eqref{ao}
to $\frac{5}{6}<\sigma < 1$, 
if the exponent $\frac{2}{3}$ will be given also for quaternion $\Gam$'s.

In the present paper, we improve and generalize \eqref{ao} as follows. 
\begin{thm}\label{thm}
Let $\Gam$ be a (not necessarily congruence) subgroup of the modular group $\sz$
or a co-compact arithmetic group $\Gamma_{\mathcal{O}}$. 
Then, for $\frac{1}{2}<\sigma<1$, we have 
\begin{align*}
\frac{1}{T}\int_{1}^{T}\left|\frac{Z'_{\Gam}(\sigma+it)}{Z_{\Gam}(\sigma+it)}\right|^{2}dt
\ll_{\epsilon,\Gam} T^{\max(0, 5-6\sigma+\epsilon)}, 
\as T\tinf.
\end{align*}
Especially, when $\frac{5}{6}<\sigma<1$, it holds
\begin{align*}
\int_{1}^{T}\left|\frac{Z'_{\Gam}(\sigma+it)}{Z_{\Gam}(\sigma+it)}\right|^{2}dt
\sim C_{\Gam}(\sigma)T, \as T\tinf
\end{align*}
for some constant $C_{\Gam}(\sigma)>0$ depending on $\Gam$ and $\sigma$.
\end{thm} 

In the theorem above, the range $\frac{5}{6}<\sigma<1$  
is available without the assumption of the Lindel\"{o}f conjecture. 
To prove the theorem, 
we first state an explicit formula (Proposition \ref{prop1}) 
to describe $Z'_{\Gam}(s)/Z_{\Gam}(s)$ 
as a sum of $\gam\in \Prim(\Gam)$ 
and take its square integral. 
In the process of estimating the square integral, 
we use an upper bound (Lemma \ref{lemclass}) of the number $m_{\Gam}(n)$ 
which is almost the same as the number of $\gam\in \Prim(\Gam)$ with the same $N(\gam)$, 
i.e. the multiplicity in the length spectrum on $\Gam\bsla H$.  
Since its multiplicity is more informative than the error term of the prime geodesic theorem, 
Theorem \ref{thm} gives an improvement of the contribution of \eqref{ao}.

\section{Explicit formula} 

To prove Theorem \ref{thm}, 
we first state the following explicit formula for the logarithmic derivative of
Selberg's zeta function. 
\begin{prop}\label{prop1} 
Let $\Gam$ be a discrete subgroup of $\sr$ with $\vol(\Gam\bsla H)<\infty$, 
$s=\sigma+iT\in \bC$ with $1/2<\sigma<1$ and 
$$
\varphi_{s}(x):=\sum_{\begin{subarray}{c}\gam\in\Prim(\Gam),j\geq1 \\ N(\gam)^j<x\end{subarray}} 
\Lambda_{\Gam}(\gam^{j})
\left(1-\frac{N(\gam)^{j}}{x} \right)N(\gam)^{-js} 
$$
for $x>0$, where $\disp \Lambda_{\Gam}(\gam^{j}):=\frac{\log{N(\gam)}}{1-N(\gam)^{-j}}$. 
Then we have 
\begin{align*}
\frac{Z'_{\Gam}(s)}{Z_{\Gam}(s)}=
&\varphi_{s}(x) 
-\sum_{\begin{subarray}{c} \frac{1}{2}<\rho\leq 1 \\ Z_{\Gam}(\rho)=0\end{subarray}}
\frac{x^{\rho-s}}{(\rho-s)(1+\rho-s)}
+O_{\epsilon}\left(T^{1+\epsilon}x^{\frac{1}{2}-\sigma+\epsilon}\right), 
\as T,x\tinf
\end{align*}
for any $\epsilon>0$, 
where $\{1/2<\rho\leq 1\}$ is the set of (at most a finite number of) real zeros of $Z_{\Gam}(s)$. 
\end{prop}

\noi{\it Proof.} 
For $V=T^{3}$ and $\epsilon>0$, let $C$ be the rectangle 
with the corners $1+\epsilon-iV$, $1+\epsilon+iV$, $1/2+\epsilon+iV$, $1/2+\epsilon-iV$, 
and 
\begin{align*}
J:=\frac{1}{2\pi i}\int_{\partial C}\frac{Z'_{\Gam}(z)}{Z_{\Gam}(z)}\frac{x^{z-s}}{(z-s)(z-s+1)}dz,
\end{align*}
where the integral is in an anti-clockwise direction.
It is known that the singular points of $Z'_{\Gam}(z)/Z_{\Gam}(z)$ in $C$ 
are a single pole at $z=1$ and 
at most a finite number of single poles on the real line. 
Then, due to the residue theorem, we have 
\begin{align}
J=\frac{Z'_{\Gam}(s)}{Z_{\Gam}(s)}+\sum_{\begin{subarray}{c} \frac{1}{2}<\rho \leq 1 \\ 
Z_{\Gam}(\rho)=0\end{subarray}}
\frac{x^{\rho-s}}{(\rho-s)(1+\rho-s)}.
\end{align}
Next, divide $J$ by 
\begin{align*}
2\pi i J=\int_{\partial C}=&\int_{1+\epsilon-i\infty}^{1+\epsilon+i\infty}
-\int_{1+\epsilon-i\infty}^{1+\epsilon-iV}-\int_{1+\epsilon+iV}^{1+\epsilon+i\infty}
+\int_{1+\epsilon+iV}^{1/2+\epsilon+iV}+\int_{1/2+\epsilon+iV}^{1/2+\epsilon-iV}
+\int_{1/2+\epsilon-iV}^{1+\epsilon-iV}\\
=:&J_{1}-J_{2}-J_{3}+J_{4}+J_{5}+J_{6}.
\end{align*}
Due to \eqref{SZorg}, we can bound $J_{2},\dots,J_{6}$ by 
\begin{align}
J_{2},J_{3}
&\ll_{\epsilon}\int_{V}^{\infty}|u|^{\epsilon}\frac{x^{1-\sigma}}{(u-T)^{2}+1}du
\ll V^{-1+\epsilon}x^{1-\sigma},\\
J_{4},J_{6}
&\ll_{\epsilon} V^{1+\epsilon}x^{1-\sigma}(V-T)^{-2}\ll V^{-1}x^{1-\sigma},\\
J_{5}&\ll_{\epsilon} \int_{-V}^{V}|u|^{1+\epsilon}\frac{x^{1/2-\sigma}}{(u-T)^{2}+1}du
\ll_{\epsilon} T^{1+\epsilon}x^{1/2-\sigma}.
\end{align}
The remaining term $J_{1}$ is given by 
\begin{align}
\frac{1}{2\pi i}J_{1}=& \sum_{\begin{subarray}{c}\gam\in\Prim(\Gam),j\geq1 \end{subarray}} 
\Lambda_{\Gam}(\gam^{j})N(\gam)^{-js}
\cdot \frac{1}{2\pi i}\int_{1-s+\epsilon-i\infty}^{1-s+\epsilon+i\infty} 
\frac{(x/N(\gam)^{j})^{z}}{z(z+1)}dz\notag\\
=&\sum_{\begin{subarray}{c}\gam\in\Prim(\Gam),j\geq1 \\ N(\gam)^j<x\end{subarray}} 
\Lambda_{\Gam}(\gam^{j})
\left(1-\frac{N(\gam)^{j}}{x} \right)N(\gam)^{-js}=\varphi_{s}(x).
\end{align} 
The proposition follows from (2.1)--(2.5) immediately. \qed

\section{Length spectrum} 

Since $N(\gam)^{j/2}=\frac{1}{2}\left(\tr{\gam^{j}}+ \sqrt{(\tr{\gam^{j}})^{2}-4} \right)$, 
the function $\varphi_{s}(x)$ 
can be expressed by 
\begin{align}\label{phi}
\varphi_{s}(x)=&
 \sum_{\begin{subarray}{c} n\in \Tr{(\Gam)} \\ n<X \end{subarray}} 
m_{\Gam}(n) \Lambda(n,x)\left(1-\frac{\epsilon(n)^{2}}{x}\right)\epsilon(n)^{-2s}
= \sum_{\begin{subarray}{c} n\in \Tr{(\Gam)} \\ n<X \end{subarray}} 
m_{\Gam}(n) \bar{\Lambda}(n,x)\epsilon(n)^{-2s}, 
\end{align}
where 
\begin{align*}
X:=&x^{1/2}+x^{-1/2}, \hspace{96pt} 
\epsilon(n):=\frac{1}{2}(n+\sqrt{n^{2}-4}), \\
\bLambda(n):=&\frac{2\log{\epsilon(n)}}{1-\epsilon(n)^{-2}}, \hspace{96pt}
\bLambda(n,x):=\frac{2\log{\epsilon(n)}}{1-\epsilon(n)^{-2}}
\left(1-\frac{\epsilon(n)^{2}}{x}\right), \\
\Tr(\Gam):=&\{\tr{\gam^{j}} \divset \gam\in\Prim(\Gam),j\geq1\}, \qquad 
m_{\Gam}(n):=\sum_{\begin{subarray}{c} \gam\in\Prim(\Gam),j\geq1 \\ \tr{\gam^{j}}=n \end{subarray}} 
\frac{1}{j}.
\end{align*}
In this section, we study $\Tr(\Gam)$ and $m_{\Gam}(n)$, 
which gives a variant of the length spectrum on $\Gam\bsla H$, 
to understand $\varphi_{s}(x)$ in detail. 
When $\Gam=\sz$, it is known that $\Tr(\Gam)=\bZ_{\geq3}$ and 
$m_{\sz}(n)$ is written in the term of 
primitive indefinite binary quadratic forms (see e.g. \cite{Ga,Sa1}), 
since there is a deep connection between $\Prim(\sz)$ and
equivalence classes of binary quadratic forms. 

Let 
$Q(x,y)=[a,b,c]:=ax^2+bxy+cy^2$
be a binary quadratic form over $\bZ$ with $a,b,c\in\bZ$ and $\gcd(a,b,c)=1$, 
and $D=D(Q):=b^2-4ac$ the discriminant of $[a,b,c]$. 
We call that two quadratic forms $Q$ and $Q'$ are equivalent ($Q\sim Q'$) 
if there exists $g\in\sz$ such that $Q(x,y)=Q'\big((x,y).g\big)$. 
Denote by $h(D)$ the number of equivalence classes of the quadratic forms of given $D=b^2-4ac$. 
It is known that, if $D>0$, then there are infinitely many positive solutions $(t,u)$ 
of the Pell equation $t^2-Du^2=4$. 
Put $(t_j,u_j)=(t_j(D),u_j(D))$ the $j$-th positive solution of $t^2-Du^2=4$ 
and $\epsilon_j(D):=\frac{1}{2}(t_j(D)+u_j(D)\sqrt{D})$. 
Note that $\epsilon_1(D)$ is called the fundamental unit of $D$ in the narrow sense, 
and it holds that $\epsilon_j(D)=\epsilon_{1}(D)^j$.
For a quadratic form $Q=[a,b,c]$ and a positive solution $(t,u)$ of $t^2-Du^2=4$, let 
\begin{align}\label{1to1}
\gam\big(Q,(t,u)\big):=\begin{pmatrix}\disp\frac{t+bu}{2}&-cu\\ au& \disp\frac{t-bu}{2}\end{pmatrix}\in\slt(\bZ).
\end{align}
Conversely, for $\gam=(\gam_{ij})_{1\leq i,j\leq 2}\in \sz$, we put
\begin{align}\label{ttod}
\begin{split}
&t_{\gam}:=\gam_{11}+\gam_{22},\quad u_{\gam}:=\gcd{(\gam_{21},\gam_{11}-\gam_{22},-\gam_{12})},\\
&a_{\gam}:=\frac{\gam_{21}}{u_{\gam}},\quad b_{\gam}:=\frac{\gam_{11}-\gam_{22}}{u_{\gam}},
\quad c_{\gam}:=-\frac{\gam_{12}}{u_{\gam}}, \\
&Q_{\gam}:=[a_{\gam},b_{\gam},c_{\gam}],\quad 
D_{\gam}:=\frac{t_{\gam}^2-4}{u_{\gam}^2}=b_{\gam}^2-4a_{\gam}c_{\gam}.
\end{split}
\end{align}
It is known that \eqref{1to1} and \eqref{ttod} give a one-to-one correspondence between
equivalence classes of primitive indefinite binary quadratic forms with $D>0$ 
and primitive hyperbolic conjugacy classes of $\sz$ 
(see \cite{Sa1} and Chap. 5 in \cite{Ga}). 
Then $m_{\sz}(n)$ is described as follows \cite{Sa1,HashLS}.
\begin{align}\label{LSsz}
m_{\sz}(n)=\sum_{u\in U(n)}\frac{1}{j_{n,u}}h\left(D_{n,u} \right),
\end{align}
where $\disp D_{n,u}:=\frac{n^2-4}{u^2}$, 
$U(n):=\{u\geq1 \divset u^2\mid n^2-4, D_{n,u}\equiv 0,1\bmod{4}\}$ and 
$j=j_{u,n}\geq 1$ is an integer with 
$\epsilon_j\left( D_{n,u}\right)=\frac{1}{2}\left(n+\sqrt{n^2-4}\right)$.

According to Arakawa-Koyama-Nakasuji's work \cite{AKN}, 
we see that $\Tr{(\Gam)}\subset \bZ_{\geq3}$ 
and $m_{\Gam}(n)$ can be expressed as follows when $\Gam$ is quaternion. 
\begin{align}\label{LSquat}
m_{\Gam}(n)=\sum_{u\in U(n)}\frac{1}{j_{n,u}}h\left(D_{n,u} \right)\lambda_{\Gam}(D_{n,u}),
\end{align}
where the constant $\lambda_{\Gam}(D)\geq0$ is bounded by 
the constant depending on $\Gam$, not on $D$
(see \cite{AKN} for $\lambda_{\Gam}(D)$ in detail). 

Furthermore, when $\Gam$ is a subgroup of $\bar{\Gam}=\sz$ or $\Gam_{\cO}$ with
$[\bar{\Gam},\Gam]<\infty$, 
we see that $\Tr{(\Gam)}\subset \Tr{(\bar{\Gam})}\subset \bZ_{\geq3}$ 
and we can obtain the following formula from 
Venkov-Zograf's formula \cite{VZ}.  
\begin{align*}
m_{\Gam}(n)=\sum_{\begin{subarray}{c}\gam\in\Prim(\bar{\Gam}),j\geq1 \\ 
\tr{\gam^{j}}=n \end{subarray}}
\frac{1}{j}\tr\left(\left(\Ind_{\Gam}^{\bar{\Gam}}1\right)(\gam^{j})\right). 
\end{align*}
Since $0\leq \tr\left(\left(\Ind_{\Gam}^{\bar{\Gam}}1\right)(\gam^{j})\right)\leq [\bar{\Gam}:\Gam]$, 
we have 
\begin{align}\label{LSsub}
0\leq m_{\Gam}(n)\leq [\bar{\Gam}:\Gam]m_{\bar{\Gam}}(n).
\end{align}

\vcpt

Due to \eqref{LSsz}, \eqref{LSquat} and \eqref{LSsub}, 
we have the following lemma.
\begin{lem}\label{lemclass}
When $\Gam$ is a subgroup of the modular group $\sz$ or a quaternion $\Gam_{\cO}$, 
we have
\begin{align}\label{LSineq}
m_{\Gam}(n)\ll_{\epsilon,\Gam} n^{1+\epsilon}, \as n\tinf.
\end{align}
\end{lem}
 
\noi{\it Proof.} 
According to \eqref{LSsz}, \eqref{LSquat} and \eqref{LSsub}, 
we see that $0\leq m_{\Gam}(n)\leq \alpha_{\Gam}m_{\sz}(n)$ holds 
for some constant $\alpha_{\Gam}$. 
Then checking \eqref{LSineq} only for $\Gam=\sz$ is enough to prove Lemma \ref{lemclass}.

For the discriminant $D>0$, 
the following class number formula holds (see e.g. \cite{Cohn}).
$$
h(D)=\frac{\sqrt{D}}{\log{\epsilon_{1}{(D)}}}L(1,\chi_{D}),
$$
where $L(z,\chi_{D}):=\sum_{n\geq1} \left(\frac{D}{n}\right)n^{-z}$ ($\Re{z}\geq1$)
is the Dirichlet $L$ function. 
It is well-known that 
$L(1,\chi_{D})\ll_{\epsilon} D^{\epsilon}$ as $D\tinf$, 
and $\epsilon_{1}(D)\geq \epsilon_{1}(5)=\frac{1}{2}(3+\sqrt{5})$.
We thus obtain  
\begin{align}
m_{\sz}(n)\ll_{\epsilon}\sum_{\begin{subarray}{c} u\geq1 \\ u^{2}\mid n^{2}-4 
 \end{subarray}}\left(\frac{n^{2}-4}{u^{2}} \right)^{1/2+\epsilon}
 \ll_{\epsilon} n^{1+\epsilon}, \as n\tinf. 
\end{align}
\qed

Remark that \eqref{LSineq} is not far from the best possible upper bound of $m_{\Gam}(n)$. 
In fact, the asymptotic formula 
\begin{align*}
\sum_{n<x}m_{\Gam}(n)^{k}\sim c_{\Gam}^{(k)}\frac{x^{k+1}}{(\log{x})^{k}}, 
\as x\tinf
\end{align*} 
has been given for some constant $c_{\Gam}^{(k)}>0$ 
when $\Gam$ is quaternion and $k=2$, 
and when $\Gam$ is a congruence subgroup of $\sz$ and any $k\geq 2$ 
\cite{BLS,Pe,Lu,HashLS}. 
This means that the growth of $m_{\Gam}(n)$ is averagely close to $\disp \frac{n}{\log{n}}$ 
and then \eqref{LSineq} works well in the proof of Theorem \ref{thm}.

\section{Proof of Theorem \ref{thm}} 

According to Proposition \ref{prop1}, we have
\begin{align}
\frac{1}{T}\int_{1}^{T}\left|\frac{Z'_{\Gam}(\sigma+it)}{Z_{\Gam}(\sigma+it)}\right|^{2}dt
=&\frac{1}{T}\int_{1}^{T} |\varphi_{\sigma+it}(x)|^{2}dt
+O(T^{-1}x^{2-2\sigma})+ O(T^{2}x^{1-2\sigma+\epsilon}). 
\end{align}
We now study the first term in the right hand side. 
\begin{align*}
\frac{1}{T}\int_{1}^{T}|\varphi_{\sigma+it}(x)|^{2}dt
= & \sum_{\begin{subarray}{c}n_{1},n_{2}\in \bZ \\ 3\leq n_{1},n_{2}<X \end{subarray}}
m(n_{1})m(n_{2})\bLambda(n_{1},x)\bLambda(n_{2},x) 
\epsilon(n_{1})^{-2\sigma}\epsilon(n_{2})^{-2\sigma}
\frac{1}{T}\left| 
\int_{1}^{T}\left(\frac{\epsilon(n_{1})}{\epsilon(n_{2})}\right)^{-2it}dt 
\right|.
\end{align*}
Divide the sum in the right hand side above by 
$\sum_{n_{1}=n_{2}}+\sum_{n_{1}\neq n_{2}}=:S_{1}+S_{2}$. 
The first sum $S_{1}$ is as follows.
\begin{align*}
S_{1}=&\sum_{\begin{subarray}{c} 3\leq n <X \end{subarray}}
m_{\Gam}(n)^{2}\bLambda(n,x)^{2} 
\epsilon(n)^{-4\sigma}\\
=& \sum_{\begin{subarray}{c} 3\leq n <X \end{subarray}}
m_{\Gam}(n)^{2}\bLambda(n)^{2} 
\epsilon(n)^{-4\sigma}
+x^{-1}\sum_{\begin{subarray}{c} 3\leq n <X \end{subarray}}
m_{\Gam}(n)^{2}\bLambda(n)^{2}
\left(\frac{\epsilon(n)^{2}}{x}-2\right) 
\epsilon(n)^{2-4\sigma}.
\end{align*}
Since $m_{\Gam}(n)\ll n^{1+\epsilon}$, 
$\epsilon(n)=n+O(n^{-1})$ and $|\bLambda(n)|\ll \log{n}$, 
we have  
\begin{align*}
S_{1}\ll_{\epsilon} \sum_{3\leq n<X}n^{2-4\sigma+\epsilon}
+x^{-1}\sum_{3\leq n<X}n^{4-4\sigma+\epsilon}
\ll_{\epsilon} X^{\max(0,3-4\sigma+\epsilon)}
\ll x^{\max(0,\frac{3}{2}-2\sigma+\epsilon)}
\end{align*}
and, especially when $\sigma>3/4$, it converges as $x\tinf$ by 
$$
\lim_{x\tinf}S_{1}=\sum_{\begin{subarray}{c} n\geq 3 \end{subarray}}
m_{\Gam}(n)^{2}\bLambda(n)^{2} 
\epsilon(n)^{-4\sigma}.
$$ 

The remaining part of the proof is the estimation of $S_{2}$. 
Since 
$$ 
\int_{1}^{T}\left(\frac{\epsilon(n_{1})}{\epsilon(n_{2})}\right)^{-2it}dt 
\ll \frac{1}{\left|\log{\frac{\epsilon(n_{1})}{\epsilon(n_{2})}}\right|}
\ll \begin{cases} 
\disp \frac{n_{1}}{n_{2}-n_{1}}, &(n_{1}<n_{2}<2n_{1}),\\
1, & (n_{2}\geq n_{1}),
\end{cases}
$$
we have 
\begin{align}
S_{2}\ll& T^{-1}\sum_{3\leq n_{1}<X}m_{\Gam}(n_{1})\bLambda(n_{1},x)\epsilon(n_{1})^{-2\sigma}
\sum_{n_{1}<n_{2}<X}m_{\Gam}(n_{2})\bLambda(n_{2},x)\epsilon(n_{2})^{-2\sigma}
\frac{1}{\left|\log{\frac{\epsilon(n_{1})}{\epsilon(n_{2})}}\right|}\notag \\
\ll_{\epsilon}& T^{-1} 
\sum_{\begin{subarray}{c}n_{1}\in \Tr(\Gam)\\ 2<n_{1}<X \end{subarray}}
n_{1}^{1-2\sigma+\epsilon}\left(
\sum_{\begin{subarray}{c}n_{2}\in \Tr(\Gam)\\ n_{1}<n_{2}<2n_{1} \end{subarray}}
n_{2}^{1-2\sigma+\epsilon}\frac{n_{1}}{n_{2}-n_{1}}+
\sum_{\begin{subarray}{c}n_{2}\in \Tr(\Gam)\\ 2n_{1}\leq n_{2}<X \end{subarray}}
n_{2}^{1-2\sigma+\epsilon}\right)\notag \\
\ll_{\epsilon}&T^{-1}X^{4-4\sigma+\epsilon}\ll T^{-1}x^{2-2\sigma+\epsilon}. \label{S2}
\end{align}

We thus obtain 
\begin{align*}
\frac{1}{T}\int_{1}^{T}\left|\frac{Z'_{\Gam}(\sigma+it)}{Z_{\Gam}(\sigma+it)}\right|^{2}dt
\ll_{\epsilon}&x^{\max(0,\frac{3}{2}-2\sigma+\epsilon)}
+T^{-1}x^{2-2\sigma}+ T^{2}x^{1-2\sigma+\epsilon}. 
\end{align*}
The main theorem immediately follows from the above with $x=T^{3}$, 
and the constant $C_{\Gam}(\sigma)$ for $\frac{5}{6}<\sigma<1$ is given by 
$$
C_{\Gam}(\sigma)=\sum_{\begin{subarray}{c} n\geq 3 \end{subarray}}
m_{\Gam}(n)^{2}\bLambda(n)^{2} 
\epsilon(n)^{-4\sigma}.
$$

\qed

\vcpt

\noi{\bf Acknowledgment.} 
The author was supported by JST CREST no.JPMJCR2113 
and JSPS Grant-in-Aid for Scientific Research (C) 
no. 17K05181.

\end{document}